\documentclass[12pt, draft]{article} 
 \usepackage{dsfont}
 \usepackage{ams}
 \usepackage{latexsym}
   \font\twelvebm                       = cmmib10 at 12truept
   \font\tenbm                          = cmmib10 at 10truept
   \font\sevenbm                        = cmmib10 at 7truept
 = \tenbm \scriptscriptfont9              = \sevenbm
 
   at 10truept  at 10truept  at
 10truept  at 10truept
  at 10truept

 \mathchardef \BGamma            = "0900 \mathchardef \BDelta
 = "0901 \mathchardef \BTheta            = "0902 \mathchardef
 \BLambda           = "0903 \mathchardef \BXi               = "0904
 \mathchardef \BPi               = "0905 \mathchardef \BSigma
 = "0906 \mathchardef \BUpsilon          = "0907 \mathchardef \BPhi
 = "0908 \mathchardef \BPsi              = "0909 \mathchardef
 \BOmega            = "090A \mathchardef \Balpha            = "090B
 \mathchardef \Bbeta             = "090C \mathchardef \Bgamma
 = "090D \mathchardef \Bdelta            = "090E \mathchardef
 \Bepsilon          = "090F \mathchardef \Bzeta             = "0910
 \mathchardef \Beta              = "0911 \mathchardef \Btheta
 = "0912 \mathchardef \Biota             = "0913 \mathchardef
 \Bkappa            = "0914 \mathchardef \Blambda           = "0915
 \mathchardef \Bmu               = "0916 \mathchardef \Bnu
 = "0917 \mathchardef \Bxi               = "0918 \mathchardef \Bpi
 = "0919 \mathchardef \Brho              = "091A \mathchardef
 \Bsigma            = "091B \mathchardef \Btau              = "091C
 \mathchardef \Bupsilon          = "091D \mathchardef \Bphi
 = "091E \mathchardef \Bchi              = "091F \mathchardef \Bpsi
 = "0920 \mathchardef \Bomega            = "0921 \mathchardef
 \Bvarepsilon       = "0922 \mathchardef \Bvartheta         = "0923
 \mathchardef \Bvarpi            = "0924 \mathchardef \Bvarrho
 = "0925 \mathchardef \Bvarsigma         = "0926 \mathchardef
 \Bvarphi           = "0927
 \mathchardef \bA        = "0941 \mathchardef \bB        = "0942
 \mathchardef \bC        = "0943 \mathchardef \bD        = "0944
 \mathchardef \bE        = "0945 \mathchardef \bF        = "0946
 \mathchardef \bG        = "0947 \mathchardef \bH        = "0948
 \mathchardef \bI        = "0949 \mathchardef \bJ        = "094A
 \mathchardef \bK        = "094B \mathchardef \bL        = "094C
 \mathchardef \bM        = "094D \mathchardef \bN        = "094E
 \mathchardef \bO        = "094F \mathchardef \bP        = "0950
 \mathchardef \bQ        = "0951 \mathchardef \bR        = "0952
 \mathchardef \bS        = "0953 \mathchardef \bT        = "0954
 \mathchardef \bU        = "0955 \mathchardef \bV        = "0956
 \mathchardef \bW        = "0957 \mathchardef \bX        = "0958
 \mathchardef \bY        = "0959 \mathchardef \bZ        = "095A
 \mathchardef \ba        = "0961 \mathchardef \bb        = "0962
 \mathchardef \bc        = "0963 \mathchardef \bd        = "0964
 \mathchardef \bee       = "0965 
 \mathchardef \bff       = "0966 \mathchardef \bg        = "0967
 \mathchardef \bh        = "0968
 \mathchardef \bj        = "096A \mathchardef \bk        = "096B
 \mathchardef \bl        = "096C \mathchardef \bm        = "096D
 \mathchardef \bn        = "096E \mathchardef \bo        = "096F
 \mathchardef \bp        = "0970 \mathchardef \bq        = "0971
 \mathchardef \br        = "0972 \mathchardef \bs        = "0973
 \mathchardef \bt        = "0974 \mathchardef \bu        = "0975
 \mathchardef \bv        = "0976 \mathchardef \bw        = "0977
 \mathchardef \bx        = "0978 \mathchardef \by        = "0979
 \mathchardef \bz        = "097A

 \font\tencb            = cmssbx10 scaled \magstep4 \font\eigcb
 = cmssbx10 scaled \magstep2 \textfont8             = \tencb
 \mathchardef\bAs       = "1841
 \def\Asem#1#2{\mathop{\vrule height10.5pt depth5.5pt width0pt\bAs}_{#1}^{#2}}
 \def\asem#1#2{
          \ifmmode
         \ifinner
            \raise0.9pt\hbox{$\scriptstyle\bAs$}_{#1}^{#2}
         \else
            \Asem{#1}{#2}
         \fi
          \fi
          }

 \newtheorem{theo}{\small\bf Theorem}
 \newtheorem{lem}{\small\bf Lemma}
 \newtheorem{rem}{\small\bf Remark}
 \newenvironment{REM}{\begin{rem} \rm}{\end{rem}}
 \newtheorem{exam}{\small\bf Example}
 \newenvironment{EXAM}{\begin{exam} \rm}{\end{exam}}
  \newtheorem{defi}{\small\bf Definition}
 
 \newtheorem{cor}{\small\bf Corollary}
 \renewcommand{\Pr}{\mbox{\rm  \hspace*{.2ex}I\hspace{-.5ex}P\hspace*{.2ex}}}

 \newcommand{\be}{\begin{equation}}
 \newcommand{\ee}{\end{equation}}
 \newcommand{\law}{\stackrel{\mbox{\footnotesize d}}{=}}

 \newcommand{\E}{\mbox{\rm \hspace*{.2ex}I\hspace{-.5ex}E\hspace*{.2ex}}}

 \newenvironment{pr}[1]{{\small\bf {#1}:}}{}

 \renewcommand{\R}{\mathds{R}}

 \title{ \Large\bf Optimal moment inequalities
 for order statistics from
 nonnegative random variables}
 \author{\normalsize\bf
 {\rm by}
 \\
 Nickos Papadatos
 }
 \date{\small Department of Mathematics,
 National and Kapodistrian
 University of Athens,  Panepistemiopolis, 157 84
 Athens,
 Greece.}

\begin{document}

 \maketitle
 \vspace*{-2em}

 \begin{abstract}
 \noindent
 We obtain the best possible upper bounds for the moments
 of a single order statistic from independent, non-negative
 random variables, in terms of the population mean.
 The main result covers the independent identically distributed case.
 Furthermore, the case of the sample minimum for merely independent
 (not necessarily identically distributed) random variables
 is treated in detail.
 \vspace{3ex}

 \end{abstract}

 \thispagestyle{empty}

 \section{Introduction}
 \setcounter{equation}{0}
 The investigation of the behavior of expectations
 of order statistics in a random sample has a long history, since
 the order statistics have several applications in statistics
 and reliability. The earliest results in this direction
 are those by Placket (1947),
 concerning the sample range, followed by the well-known
 papers by Hartley and David (1954) and Gumbel (1954), regarding
 the expected extremes. At those years, a pioneer paper
 by Moriguti (1953) established a powerful projection method,
 making possible to evaluate tight expectation  bounds
 for the non-extreme order statistics in terms of the population mean
 and variance. Since then, a large number of generalizations extensions
 and improvements have been found, including linear estimators
 from dependent samples
 (Arnold and Groeneveld (1979); Rychlik 1992, 1993a, 1993b, 1998;
 Balakrishnan 1990;
 Gascuel and Caraux 1992; Papadatos 2001a; Papadatos and Rychlik 2004;
 Miziula and Navarro 2018),
 record values and $k$th records (Raqab 2004; Raqab and Rychlik 2002)
 as well as distribution bounds
 (Caraux and Gascuel 1992; Papadatos 2001b, Okolewski 2015), to mention
 a few.
 The reader is referred to the monographs by Arnold and Balakrishnan
 (1989), Rychlik (2001) and Ahsanullah and Raqab (2006) for
 a comprehensive presentation on characterizations and bounds through order
 statistics and records.

 Beyond the well-developed theory on expectation bounds for order
 statistics and records, the corresponding theory
 to other moments does not seem to have receive much attention.
 Of course, some exceptions exist concerning variances; see, e.g.,
 Papadatos (1995), Jasi\'nski and Rychlik (2012, 2016),
 Rychlik (2008, 2014). The purpose of the present work is to obtain
 tight upper bounds for the moments of a single order statistic
 from a nonnegative population. These bounds are useful at least
 for reliability systems, since, as is well-known, the $k$th order
 statistic, $X_{k:n}$, represents the time-to-failure
 in a $(n+1-k)$-out-of-$n$ system -- clearly, the individual
 components cannot have negative lifetimes, hence the assumption of
 nonnegativity is natural
 for this kind of problems.

 The paper is organized as follows. In section 2 we provide general results
 for the existence of moments of a single order statistic
 in the general (not necessarily identically distributed) independent case.
 Section 3 contains the main results,
 providing tight upper bounds for the moments of
 non-extreme order statistics in terms of the population mean.
 Finally, Section 4 contains detailed results for the sample minimum,
 which represents the lifetime of a serial system.

 \section{Moment bounds in the independent case}
 \setcounter{equation}{0}

 Let $X_1,\ldots,X_n$ be $n$ iid (independent identically distributed)
 copies of the random variable (rv)
 $X$ and consider
 the corresponding order statistics $X_{1:n}\leq \cdots\leq X_{n:n}$.
 It is well-known that if $X$ is integrable then the same is
 true for any order statistic $X_{i:n}$ (for all $n$ and $i$).
 Moreover,
 an old result by P.K.\ Sen (1959) showed that the condition
 \[
 \E|X|^\delta <\infty \ \ \mbox{ for some } \ \delta\in(0,1]
 \]
 is sufficient for
 \[
 \E|X_{i:n}|<\infty \ \ \mbox{ for all $i$ with } \ \frac{1}{\delta}\leq i\leq
 n+1-\frac{1}{\delta}.
 \]

 It is natural to look at similar conditions when $X$ is nonnegative
 (cf.\ Papadatos, 1997), since this is the case for several applications
 including reliability systems. The main purpose of the present work
 is in obtaining best possible bounds for the moments of a single
 order statistic
 from non-negative populations, in terms of the population mean.

 The results of the present section concern the existence of moments
 in the more general case where the $X_i$'s are merely independent.
 We have the following.

 \begin{lem}
 \label{lem.1}
 If $X_1,\ldots,X_n$ are non-negative independent rv's
 with $\E X_i<\infty$ then $\E (X_{k:n})^{n+1-k}<\infty$. In particular,
 $X_{1:n}$ has finite $n$-th moment, $X_{2:n}$ has finite $(n-1)$-th moment,
 and $X_{n:n}$ has finite first moment.
 \end{lem}

 \noindent
 The proof of Lemma \ref{lem.1} is evident from Theorem \ref{theo.1}
 below.

 \begin{lem}
 \label{lem.2}
 Given $\mu_1,\ldots,\mu_n>0$,
 there are non-negative independent rv's $X_1,\ldots,X_n$ with
 $\E X_i=\mu_i$ for all $i$ and $\E (X_{k:n})^{n+1-k+\delta}=\infty$
 for all $k\in\{1,\ldots,n\}$ and for any $\delta>0$.
 Moreover, if $\mu_1=\cdots=\mu_n$, the
 rv's $X_1,\ldots,X_n$ can be chosen to be iid.
 \end{lem}

 \noindent
 \begin{pr}{Proof}
 For $\mu>0$ consider the function
 \[
 R_{\mu}(x)=\left\{
 \begin{array}{ll}
 1, & \mbox{ if } x\leq \mu/2,
 \\
 \displaystyle
 \frac{1}{\frac{2x}{\mu}(1+\log \frac{2x}{\mu})^2},
 & \mbox{ if } x\geq \mu/2.
 \end{array}
 \right.
 \]
 It is easy to check that $R_\mu(x)$ is a reliability function
 of an rv, $Y_\mu$, say, that is, $Y_\mu\sim F_{\mu}=1-R_{\mu}$.
 Obviously, $Y_{\mu}$ is supported in $(\frac{\mu}{2},\infty)$
 and, moreover, $\lambda Y_\mu\sim F_{\lambda \mu}$, $\lambda>0$;
 hence, $Y_{\mu}\law \mu Y_1$, where $\law$ denotes equality in
 distribution. Furthermore,
 \[
 \E Y_1 =\int_{0}^\infty R_1(x) dx
 =\frac{1}{2}
 +\int_{1/2}^\infty \frac{1}{2x(1+\log 2x)^2} dx
 =\frac{1}{2}
 +\frac{1}{2}\int_{0}^\infty \frac{1}{(1+t)^2} dt=1,
 \]
 where we made use of the substitution $\log 2x =t$.
 For any $\alpha>0$, a similar calculation yields
 \[
 \E (Y_1) ^\alpha
 =\alpha \int_{0}^\infty
 x^{\alpha-1} R_1(x) dx=
 \frac{1}{2^{\alpha}}
 +\frac{\alpha}{2^\alpha}\int_{0}^\infty \frac{e^{-(1-\alpha) t}}{(1+t)^2} dt;
 \]
 note that this formula holds even if $(Y_1)^\alpha$
 is non-integrable (see, e.g., Jones and Balakrishnan 2002).
 Hence, $\E (Y_1)^\alpha <\infty$ if and only if $\alpha\in(0,1]$.

 Without loss of generality assume that $0<\mu_1\leq \cdots\leq \mu_n$ and
 consider the independent rv's $X_1,\ldots,X_n$ with $X_i\law \mu_i Y_1$,
 $i=1,\ldots,n$. It is clear that the $X_i$'s are iid if and only if
 the $\mu_i$'s
 are all equal. Moreover, consider the iid rv's $Z_1,\ldots,Z_n$
 with $Z_{i}=\frac{\mu_1}{\mu_i} X_i$, $i=1,\ldots,n$. Since the function
 $(x_1,\ldots,x_n)\mapsto x_{k:n}(x_1,\ldots,x_n)$ is non-decreasing in its
 arguments and $Z_i\leq X_i$,
 we have
 \[
 Z_{k:n}\leq X_{k:n}, \ \ k=1,\ldots,n.
 \]
 Hence, it suffices to show that $\E (Z_{k:n}) ^{n+1-k+\delta}=\infty$
 for $\delta>0$.
 To this end, observe that the $Z_i$'s are iid from $F_{\mu_1}$
 and
 \[
 \Pr(Z_{k:n}>x)=\sum_{j=n+1-k}^{n} {n \choose j} R_{\mu_1}(x)^j F_{\mu_1}(x)^{n-j}
 \geq {n \choose k-1} R_{\mu_1}(x)^{n+1-k} F_{\mu_1}(x)^{k-1};
 \]
 the lower bound is just the first term of the sum. Therefore,
 since $F_{\mu_1}(x)> \frac{1}{2}$ for $x\geq \mu_1$, we have
 \begin{eqnarray*}
 \E (Z_{k:n})^{n+1-k+\delta}
 & = &
 (n+1-k+\delta)\int_{0}^{\infty}
  x^{n-k+\delta} \Pr(Z_{k:n}>x) dx
 \\
  &
 \geq
 &
 \frac{n+1-k+\delta}{2^{k-1}} {n \choose k-1}
 \int_{\mu_1}^{\infty} x^{n-k+\delta} R_{\mu_1}(x)^{n+1-k} dx
 \\
 &
 =
 &
 \frac{(n+1-k+\delta)(\mu_1)^{n+1-k+\delta}}{2^{n+\delta}} {n \choose k-1}
 \int_{\log \mu_1}^{\infty} \frac{e^{\delta t}}{(1+t)^{2(n+1-k)}} dt
 \\
 &
 =
 &
 \infty,
 \end{eqnarray*}
 completing the proof.
 $\Box$
 \bigskip
 \end{pr}

 Our result for the independent case is the following
 theorem, which includes Lemma \ref{lem.1} as a particular case.

 \begin{theo}
 \label{theo.1}
 If $X_1,\ldots,X_n$ are non-negative independent rv's
 with $\E X_i=\mu_i\in(0,\infty)$ then
 \be
 \label{1}
 \E (X_{k:n})^{n+1-k}\leq \sum_{1\leq i_1<\cdots<i_{n+1-k}\leq n} \mu_{i_1}
 \cdots \mu_{i_{n+1-k}}, \ \ k=1,\ldots,n.
 \ee
 The equality in {\rm (\ref{1})} is attainable for $k=1$, and is best possible
 and non-attainable for $k\geq 2$.
 \end{theo}
 \begin{pr}{Proof}
 Observe that
 \begin{eqnarray*}
 (X_{k:n})^{n+1-k}
 &
 \leq
 &
 X_{k:n} \cdots X_{n:n}
 \\
 &&
 \\
 &
 \leq
 &
 \sum_{1\leq i_1<\cdots<i_{n+1-k}\leq n}
 X_{i_1:n} \cdots X_{i_{n+1-k}:n}
 \\
 &
 =
 &
 \sum_{1\leq i_1<\cdots<i_{n+1-k}\leq n} X_{i_1} \cdots X_{i_{n+1-k}}.
 \end{eqnarray*}
 Hence, taking expectations and using the fact that the $X_i$'s are independent,
 we deduce (\ref{1}). We shall now verify that for $k\geq 2$ the equality
 is non-attainable. Indeed, if $k\geq 2$, the above sum
 contains at least two summands. Let
 \[
 Y_1=(X_{k:n})^{n+1-k}, \ \ \ \
 Y_2=X_{k:n} \cdots X_{n:n}, \ \  \ \ Y_3=\sum_{1\leq i_1<\cdots<i_{n+1-k}\leq n}
 X_{i_1:n} \cdots X_{i_{n+1-k}:n},
 \]
 so that $Y_1\leq Y_2\leq Y_3$. Assuming equality in (\ref{1}), that is,
 $\E Y_1=\E Y_3$, we see that $\E(Y_3-Y_2)=0$. Therefore, taking expectations
 to
 the obvious
 inequalities
 $0\leq (X_{1:n})^{n+1-k}\leq X_{1:n}\cdots X_{n+1-k:n}\leq Y_3-Y_2$ (the last one
 is valid because $k\geq 2$), we obtain
 \[
 0\leq \E (X_{1:n})^{n+1-k} \leq \E (Y_3-Y_2)=0.
 \]
 Hence, $X_{1:n}=0$ w.p.\ $1$. However, this fact is impossible, since
 \[
 \Pr(X_{1:n}>0)=\Pr(X_1>0,\ldots,X_n>0)=\prod_{j=1}^n \Pr(X_j>0)>0,
 \]
 because $\E X_i>0$. The attainability of (\ref{1}) for $k=1$ will be shown
 later (see Theorem \ref{th2}, below).

 We now show that inequality (\ref{1}) is best possible.
 Fix $M\geq \max_i\{\mu_i\}$ and
 consider independent two-valued rv's $X_i$ with
 \[
 \Pr(X_i=0)=1-\frac{\mu_i}{M},  \ \ \Pr(X_i=M)=\frac{\mu_i}{M},
 \]
 so that $\E X_i =\mu_i$ for all $i$. It is easy to see that
 \begin{eqnarray*}
 \Pr(X_{k:n}=M)
 &=&
 \Pr(\mbox{at least $n+1-k$ among $X_1,\ldots,X_n$ are equal to $M$})
 \\
 &\geq&
 \Pr(\mbox{exactly $n+1-k$ among $X_1,\ldots,X_n$ are equal to $M$})
 \\
 &=&
 \sum_{1\leq i_1<\cdots<i_{n+1-k}\leq n}
 \frac{\mu_{i_1}\cdots\mu_{i_{n+1-k}}}{M^{n+1-k}}
 \prod_{j\in S(i_1,\ldots,i_{n+1-k})} \left(1-\frac{\mu_j}{M}\right),
 \end{eqnarray*}
 where $S(i_1,\ldots,i_{n+1-k})=\{1,\ldots,n\}\setminus\{i_1,\ldots,i_{n+1-k}\}$.
 The smallest term in the product is
 at least $1-\frac{\max_{i}\{\mu_i\}}{M}$, hence,
 \[
 \Pr(X_{k:n}=M) \geq \left(1-\frac{\max_{i}\{\mu_i\}}{M}\right)^{k-1}
 \frac{1}{M^{n+1-k}}  \sum_{1\leq i_1<\cdots<i_{n+1-k}\leq n}
 \mu_{i_1}\cdots\mu_{i_{n+1-k}}.
 \]
 It follows that
 \begin{eqnarray*}
 \E (X_{k:n})^{n+1-k}
 & = &
 M ^{n+1-k} \ \Pr(X_{k:n}=M)
 \\
 &\geq&
 \left(1-\frac{\max_{i}\{\mu_i\}}{M}\right)^{k-1}
  \sum_{1\leq i_1<\cdots<i_{n+1-k}\leq n}
 \mu_{i_1}\cdots\mu_{i_{n+1-k}}
 \\
 &\to&
   \sum_{1\leq i_1<\cdots<i_{n+1-k}\leq n}
 \mu_{i_1}\cdots\mu_{i_{n+1-k}},  \ \
 \mbox{as $M\to\infty$},
\end{eqnarray*}
 and the proof is complete.
 $\Box$
 \bigskip
 \end{pr}
 \begin{cor}
 \label{cor.1}
 If $\mu_1=\cdots=\mu_n=\mu>0$ {\rm (in particular, if the $X_i$'s are iid)},
 the best possible upper bound is
 given by
 \[
  \E (X_{k:n})^{n+1-k} \leq {n\choose k-1} \mu^{n+1-k},  \ \ k=1,\ldots,n,
 \]
 and it is attainable only in the case $k=1$.
 \end{cor}

 \section{Moment bounds for the independent, identically distributed, case}

 In this section we assume that $X_1,\ldots,X_n$ are iid
 non-negative rv's distributed like $X$, and $\E X=\mu$
 is nonzero and finite. Our purpose is to derive
 the best possible upper bounds for the moments
 $\E (X_{k:n})^{\alpha}$, for $\alpha>0$;
 however, due to Lemmas \ref{lem.1} and \ref{lem.2}, we see that
 the problem is meaningful only for $\alpha\in(0,n+1-k]$.
 Note that Papadatos (1997) treats the case $\alpha=1$, which, as
 we shall see below, is a boundary case between $\alpha<1$
 and $\alpha>1$. Also, we shall obtain the populations
 that attain the equality in the bounds.

 We first prove some auxiliary results. In the following lemma
 we consider
 the usual Borel space
 \[
 L^1(0,1)=\left\{ g:(0,1)\to\R, \ g \mbox{ Borel, }
 \int_0^1 |g(t)|dt<\infty\right\},
 \]
 where two functions that differ at a set of Lebesgue
 measure zero are considered as equal.
 \begin{lem}
 \label{lem.3} Let $\alpha>1$.
 If a function $g:(0,1)\to [0,\infty)$ is nondecreasing
 and belongs to $L^1(0,1)$, then
 \be
 \label{2}
 \alpha \int_{0}^1 (1-t)^{\alpha-1} g(t)^\alpha dt \leq
 \left(\int_0^1 g(t) dt\right)^\alpha,
 \ee
 and the equality holds if
 either $g$ is constant or
 \[
 g(t) =\left\{
 \begin{array}{ll}
 0, & 0<t\leq t_0,
 \\
 \theta, & t_0< t<1,
 \end{array}
 \right.
 \]
 for some $t_0\in(0,1)$ and some $\theta>0$.
 \end{lem}
 \begin{pr}{Proof}
 It is obvious that the any constant function attains the equality in
 (\ref{2}), and the same is true for the function
 $g(t)=\theta I_{(t_0,1)}(t)$, resulting to the identity
 $\theta^\alpha (1-t_0)^\alpha=(\theta(1-t_0))^\alpha$.
 To prove the inequality,
 assume first that $g$ is simple nonnegative and nondecreasing,
 that is,
 \[
 g(t)=\left\{
 \begin{array}{ll}
 \delta_1, & t\in(0,s_1],
 \\
 \delta_1+\delta_2, & t\in(s_1,s_2],
 \\
 \vdots &
 \\
 \delta_1+\delta_2+\cdots+\delta_k, & t\in(s_{k-1},1),
 \end{array}
 \right.
 \]
 where $\delta_i\geq 0$ and $0<s_1<\cdots<s_{k-1}<1$; note that
 the value of $g$ at the end-points do not affect the value
 of the integrals, so we have assumed that $g$ is left-continuous.
 With the notation $s_0=0$, $s_k=1$,
 it is easily seen that
 \[
 \int_{0}^1 g(t) dt
 = \sum_{j=1}^k (s_j-s_{j-1}) (\delta_1+\cdots+\delta_j)
 =\sum_{j=1}^k (1-s_{j-1})\delta_j.
 \]
 Similarly,
 \[
 \alpha \int_{0}^1 (1-t)^{\alpha-1} g(t)^\alpha dt
 =
 (\delta_1)^\alpha +\sum_{j=1}^{k-1} (1-s_{j-1})^\alpha
 \left[
 (\delta_1+\cdots+\delta_{j+1})^\alpha-
 (\delta_1+\cdots+\delta_{j})^\alpha
 \right].
 \]
 Therefore, (\ref{2}) for simple functions
 reduces to the inequality
 \be
 \label{3}
 \left(\sum_{j=1}^k (1-s_{j-1})\delta_j\right)^\alpha-
 (\delta_1)^\alpha -\sum_{j=1}^{k-1} (1-s_{j-1})^\alpha
 \left[
 (\delta_1+\cdots+\delta_{j+1})^\alpha-
 (\delta_1+\cdots+\delta_{j})^\alpha
 \right]\geq 0
 \ee
 for $k\geq 2$, $0=s_0<s_1<\cdots<s_{k-1}<s_k=1$ and $\delta_j\geq 0$
 ($j=1,\ldots,k$); note that $k=1$ leads to the constant function
 $g\equiv \delta_1$, and in this case we have equality in (\ref{2}).
 We shall show (\ref{3}) using induction on $k$. For $k=2$, (\ref{3}) reads
 as
 \[
 f(\delta_2):=[\delta_1+(1-s_1)\delta_2]^\alpha - (\delta_1)^\alpha-
 (1-s_{1})^\alpha
 \left[
 (\delta_1+\delta_{2})^\alpha-
 (\delta_1)^\alpha
 \right]\geq 0.
 \]
 However, this follows easily because $f(0)=0$ and
 \[
 f'(\delta_2)=\alpha (1-s_1)  \left[
 (\delta_1+(1-s_1)\delta_2)^{\alpha-1}
 -((1-s_1) \delta_1+(1-s_1)\delta_2)^{\alpha-1}\right]\geq 0,
 \]
 since $\alpha>1$ and $\delta_1+(1-s_1)\delta_2\geq (1-s_1)
 \delta_1+(1-s_1)\delta_2$. Assuming that (\ref{3}) holds for some
 $k\geq 2$, we shall verify it for $k+1$. Set
 \begin{eqnarray*}
 f(\delta_{k+1})
 &
 :=
 &
 \left(\sum_{j=1}^{k+1} (1-s_{j-1})\delta_j\right)^\alpha-
 (\delta_1)^\alpha
 \\
 &
 &
 -\sum_{j=1}^{k} (1-s_{j-1})^\alpha
 \left[
 (\delta_1+\cdots+\delta_{j+1})^\alpha-
 (\delta_1+\cdots+\delta_{j})^\alpha
 \right].
 \end{eqnarray*}
 It is easily seen that $f(0)\geq 0$, due to the induction argument.
 Moreover,
 \[
 f'(\delta_{k+1})=\alpha(1-s_k)\left[
 \left(\sum_{j=1}^{k+1}(1-s_{j-1})\delta_j\right)^{\alpha-1}
 -\left(\sum_{j=1}^{k+1}(1-s_{k})\delta_j\right)^{\alpha-1}\right]\geq 0,
 \]
 since $\alpha>1$ and $\sum_{j=1}^{k+1}(1-s_{j-1})\delta_j
 \geq \sum_{j=1}^{k+1}(1-s_{k})\delta_j$.
 Hence, (\ref{3}) is valid for simple functions.
 If $g\geq 0$ is an arbitrary nondecreasing right-continuous
 function, we can use standard
 arguments to find simple
 functions $g_n$ such that $g_n\nearrow g$ pointwise. Then,
 $\alpha (1-t)^{\alpha-1} g_n(t)\nearrow \alpha (1-t)^{\alpha-1} g(t)$ and,
 by Lebesgue's monotone convergence theorem and (\ref{3}) we get
 \begin{eqnarray*}
 &&\alpha \int_{0}^1 (1-t)^{\alpha-1} g(t)^\alpha dt
 =\lim_n \left\{\alpha \int_{0}^1 (1-t)^{\alpha-1} g_n(t)^\alpha dt
 \right\}
 \\
 &&
 \leq
 \lim_n
 \left(\int_0^1 g_n(t) dt\right)^\alpha
 =\left(\lim_n \int_0^1 g_n(t) dt\right)^\alpha
 =\left(\int_0^1 g(t) dt\right)^\alpha,
 \end{eqnarray*}
 completing the proof.
 $\Box$
 \bigskip
 \end{pr}

 \begin{cor}
 \label{cor.2}
 Let $F$ be a distribution function of a nonnegative rv $X$
 with mean $\mu\in(0,\infty)$. Then, for all $\alpha>1$,
 \[
 \alpha \int_{0}^\infty x^{\alpha-1} (1-F(x))^\alpha dx\leq
 \left(\int_{0}^\infty (1-F(x)) dx\right)^\alpha,
 \]
 and the equality is attained if $X$ assumes two values, one
 of which is zero.
 \end{cor}
 \begin{pr}{Proof}
 It is trivial to check that any distribution function (df)
 $F(x)$
 that is constant in $[0,x_0)$ and equals to one in $[x_0,\infty)$
 attains the equality. We now verify the
 inequality;
 note that for integral values of $\alpha>1$, say $\alpha=n$,
 it becomes obvious if we consider
 the rv $X_{1:n}=\min\{X_1,\ldots,X_n\}$,
 where $X_1,\ldots,X_n$ are iid from $F$. Then,
 \[
 \E (X_{1:n})^n \leq \E (X_{1:n} \cdots X_{n:n}) =
 \E (X_{1} \cdots X_{n}) = \mu^n,
 \]
 and this inequality is equivalent to the desired one for $\alpha=n$.
 However, this simple argument is not sufficient to prove the result
 for non-integral
 values of $\alpha>1$.
 In order to verify the inequality in its general form,
 let $F^{-1}(u)=\inf\{x:F(x)\geq u\}$, $0<u<1$, be the left-continuous inverse
 of $F$. Moreover, consider an rv $Y_\alpha$ with
 df $F_{\alpha}=1-(1-F(x))^\alpha$.
 It is easy to see that $F_\alpha ^{-1}(u)=F^{-1}(1-(1-u)^{1/\alpha})$.
 Hence, from Lemma \ref{lem.3}
 with $g=F^{-1}$,
 \[
 \E (Y_{\alpha})^\alpha =\int_{0}^1 (F_{\alpha}^{-1}(u))^\alpha du
 =\alpha \int_{0}^1  (1-t)^{\alpha-1} (F^{-1}(t))^\alpha dt
 \leq \left(\int_{0}^1 F^{-1}(t) dt \right)^\alpha =\mu^\alpha,
 \]
 where we used the substitution $t=1-(1-u)^{1/\alpha}$.
 Moreover, since
 \[
 \E (Y_{\alpha})^\alpha= \alpha \int_{0}^{\infty} x^{\alpha-1} (1-F(x))^\alpha dx
 \mbox{ \ \ \ and \ \ \ }
 \mu^\alpha = \left(\int_{0}^{\infty}  (1-F(x)) dx\right)^{\alpha},
 \]
 the result if proved.
 $\Box$
 \bigskip
 \end{pr}

 \begin{lem}
 \label{lem.4}
 Let $n\geq 3$, $k\in\{2,\ldots,n-1\}$ and $\alpha\in [1,n+1-k)$.
 Let, also,
 \be
 \label{1a}
 G_{k:n}(x)=\sum_{j=k}^n {n\choose j} x^j (1-x)^{n-j}, \ \ 0\leq x\leq 1,
 \ee
 be the df of $U_{k:n}$ from an iid sample
 $U_1,\ldots,U_n$ from the standard uniform df, and
 \be
 \label{2a}
 g_{k:n}(x)=G_{k:n}'(x)=\frac{1}{B(k,n+1-k)} x^{k-1}(1-x)^{n-k}, \ \ 0<x<1,
 \ee
 the corresponding Beta density.
 Then,
 \be
 \label{4}
 1-G_{k:n}(x)\leq A_{k:n}(\alpha) (1-x)^{\alpha}, \ \ 0\leq x\leq 1,
 \ee
 where
 \be
 \label{5}
 A_{k:n}(\alpha)=\frac{1-G_{k:n}(\rho)}{(1-\rho)^\alpha}
 \ee
 and $\rho=\rho_{k:n}(\alpha)$ is the unique solution to the equation
 \be
 \label{6}
 \alpha (1-G_{k:n}(\rho))=(1-\rho) g_{k:n}(\rho), \ \ 0<\rho<1.
 \ee
 The equality in {\rm (\ref{4})} is attained if and only if $x=\rho$ or $x=1$.
 \end{lem}
 \begin{pr}{Proof}
 Define the function
 \[
 h(x)=\frac{1-G_{k:n}(x)}{(1-x)^\alpha}, \ \ 0\leq x\leq 1.
 \]
 where the value at $x=1$ is defined by continuity: $h(1)=0$.
 We have $h(0)=1$, $h(1)=0$ and
 \[
 h'(x)=(1-x)^{-\alpha-1}\Big(\alpha (1-G_{k:n}(x))-(1-x) g_{k:n}(x)\Big), \ \
 0<x<1.
 \]
 Setting $t(x)=\alpha (1-G_{k:n}(x))-(1-x) g_{k:n}(x)$,
 we calculate
 \[
 t'(x)= \frac{g_{k:n}(x)}{x} \Big((n-\alpha)x-(k-1)\Big).
 \]
 This shows that $t(x)$ is strictly decreasing in $(0,\frac{k-1}{n-\alpha}]$
 and strictly increasing in $[\frac{k-1}{n-\alpha},1)$. Since
 $t(0)>0$ and $t(1)=0$, the function $t$ has a global negative minimum
 at $\frac{k-1}{n-\alpha}$ and, therefore, there exists a
 $\rho\in(0,\frac{k-1}{n-\alpha})$ such that $t(x)>0$ for $x\in(0,\rho)$
 and $t(x)<0$ for $x\in(\rho,1)$. Since $h'(x)=t(x)/(1-x)^{\alpha+1}$,
 we see that the function $h$ is strictly increasing in $(0,\rho)$
 and strictly decreasing in $(\rho,1)$, attaining its global maximum
 at $x=\rho$, where $\rho$ is the unique root of (\ref{6}).
 $\Box$
 \bigskip
 \end{pr}

 \begin{REM}
 \label{rem.1}
 Due to (\ref{6}), we can write $\displaystyle A_{k:n}(\alpha)
 =\frac{g_{k:n}(\rho)}{\alpha(1-\rho)^{\alpha-1}}$.
 \end{REM}

 We can now state and prove
 the main result for the moments of the non-extreme
 order statistics.
 \begin{theo}
 \label{theo.2}
 Let $X_1,\ldots,X_n$ {\rm ($n\geq 3$)} be iid nonnegative rv's with
 mean $\mu\in(0,\infty)$.
 Then, for any $k\in\{2,\ldots,n-1\}$ and $\alpha\in[1,n+1-k)$,
 \be
 \label{7}
 \E (X_{k:n})^{\alpha}\leq A_{k:n}(\alpha) \ \mu ^\alpha,
 \ee
 where $A_{k:n}(\alpha)$ is given by {\rm (\ref{5})}.
 The equality in {\rm (\ref{7})} is attained if and only if
 $\Pr(X_i=0)=\rho$, $\Pr(X_i=\mu/(1-\rho))=1-\rho$, where
 $\rho=\rho_{k:n}(\alpha)$ is given by {\rm (\ref{6})}.
 \end{theo}
 \begin{pr}{Proof}
 If $F$ is the df of the $X_i$'s then
 $G_{k:n}\circ
  F$ is the df of $X_{k:n}$;
 see David (1981), Arnold {\it et al} (2008),  David and Nagaraja (2003).
 Therefore,
 \begin{eqnarray*}
 \E (X_{k:n})^{\alpha}
 & = &
 \alpha \int_{0}^{\infty} x^{\alpha-1}
 (1-G_{k:n}(F(x))) dx
 \\
 & \leq &
 A_{k:n}(\alpha) \alpha \int_{0}^{\infty} x^{\alpha-1}
 (1-F(x))^\alpha dx
 \\
 & \leq &
 A_{k:n}(\alpha) \left(\int_{0}^{\infty} (1-F(x)) dx\right)^\alpha
 \ = \ A_{k:n}(\alpha)\  \mu^\alpha,
 \end{eqnarray*}
 where the first inequality follows from Lemma \ref{lem.4}
 and the second one from Corollary \ref{cor.2}.
 In order to have equality in (\ref{7}), it is necessary and sufficient
 that the set $\{F(x), \ 0\leq x<\infty\}$ coincides with
 $\{\rho,1\}$ -- see (\ref{6}) and Corollary \ref{cor.2}.
 Therefore, $X_1$ assumes the value $0$ w.p.\ $\rho$ and
 a positive value $x_0$ w.p.\ $1-\rho$. Finally, the condition
 $\E X_1=\mu$ shows that $x_0=\mu/(1-\rho)$, completing the proof.
 $\Box$
 \bigskip
 \end{pr}

 \begin{REM}
 For $\alpha=1$, the bounds coincide with the upper bounds
 given in Papadatos (1997), Theorem 2.1.
 \end{REM}

 \begin{EXAM}
 \label{exam.1}
 For $k=2$ one finds
 \[
 \rho=\frac{\alpha}{(n-1)(n-\alpha)} \ \ \mbox{ and }
 \ \
 A_{2:n}(\alpha)=\left(1+\frac{\alpha}{n-\alpha}\right)^{n-\alpha}
 \left(1-\frac{\alpha}{n-1}\right)^{n-1-\alpha},
 \]
 $1\leq \alpha <n-1$. It is easy to verify that
 \[
 A_{2:n}(\alpha)=1+\frac{\alpha^2}{2 n^2}+o(n^{-2}) \ \
 \mbox{as} \ \  n\to\infty.
 \]
 Closed forms can be found for $k=3$ too; then (\ref{6})
 is reduced to a second degree polynomial equation -- see
 Balakrishnan (1993).
 \end{EXAM}

 We now turn to the case $\alpha<1$, showing the following result.
 \begin{theo}
 \label{theo.3}
 Let $X_1,\ldots,X_n$ be iid nonnegative rv's with mean $\mu\in(0,\infty)$.
 If $n\geq 2$ and $k\in \{2,\ldots,n\}$, then
 \be
 \label{8}
 \E (X_{k:n})^\alpha \leq A_{k:n}(\alpha) \ \mu^\alpha,  \ \
 0<\alpha< 1,
 \ee
 where
 \be
 \label{9}
 A_{k:n}(\alpha)=\left(\int_{0}^1
 \overline{g}_{k:n}(u)^\frac{1}{1-\alpha}du \right)^{1-\alpha}
 \ee
 and $\overline{g}_{k:n}$ is the derivative of the greatest convex minorant,
 $\overline{G}_{k:n}$,
 of the function $G_{k:n}$ given in {\rm (\ref{1a})}. Specifically,
 for $k\in\{2,\ldots,n-1\}$,
 \[
 \overline{g}_{k:n}(u)=g_{k:n}(\min\{u,\rho\}), \ \ 0<u<1,
 \]
 where $\rho=\rho_{k:n}$ is the unique root to the equation
 \[
 1-G_{k:n}(\rho)=(1-\rho) g_{k:n}(\rho), \ \ 0<\rho< 1,
 \]
 while $\overline{g}_{n:n}(u)=g_{n:n}(u)=n u^{n-1}$.
 The equality in {\rm (\ref{8})} is attained if and only
 if the inverse df of $X_1$ is given by
 \be
 \label{10}
 F^{-1}(u)= \mu \ \overline{g}_{k:n}(u)^\frac{1}{1-\alpha}\Big /
 \int_0^1
 \overline{g}_{k:n}(t)^\frac{1}{1-\alpha} dt, \ \ 0<u<1.
 \ee
 \end{theo}
 \begin{pr}{Proof}
 We shall apply a slight variation of the pioneer projection
 method due to
 Moriguti (1953).
 Since $X_{k:n}\law F^{-1}(U_{k:n})$ where $U_{k:n}$
 is the $k$-th order statistic from the standard uniform
 df,
 we have
 (see Moriguti 1953;
 Rychlik 2001; Ahsanullah and Raqab 2006, Lemma 3.1.1)
 \[
 \E (X_{k:n})^\alpha = \int_{0}^1 g_{k:n}(u) F^{-1}(u)^{\alpha} du
 \leq \int_{0}^1 \overline{g}_{k:n}(u) F^{-1}(u)^{\alpha} du,
 \]
 by Moriguti's inequality (the function $(F^{-1})^\alpha$ is, clearly,
 non-decreasing).
 Applying H\"older's inequality,
 \[
 \mbox{$\int f g\leq \left(\int f^p \right)^{1/p} \left(\int g^q \right)^{1/q}$}
 \ \ \ (p, q>1, \ \ 1/p+1/q=1),
 \]
 to the last integral,
 with $f=\overline{g}_{k:n}$,
 $g=(F^{-1})^\alpha$, $p=1/(1-\alpha)>1$ and $q=1/\alpha>1$,
 we obtain the inequality
 \[
 \int_{0}^1 \overline{g}_{k:n}(u) F^{-1}(u)^{\alpha} du
 \leq \left(\int_{0}^1
 \overline{g}_{k:n}(u)^\frac{1}{1-\alpha}du \right)^{1-\alpha}
 \left(\int_{0}^1 F^{-1}(u) du \right)^\alpha,
 \]
 which verifies (\ref{8}). We now examine the case
 of equality: it is well-known that for the H\"older inequality
 to hold as equality it is necessary and sufficient
 that $g^q=c\ f^p$ for some $c\geq 0$ (note that $f,g\geq 0$ in our case);
 that is, $F^{-1}(u)=c \ \overline{g}_{k:n}(u)^\frac{1}{1-\alpha}$.
 Taking into account the condition $\E X_1=\mu$ we get
 \[
 \mu=\int_{0}^1 F^{-1}(t) dt=c
 \int_{0}^1 \overline{g}_{k:n}(t)^\frac{1}{1-\alpha} dt.
 \]
 Therefore, $c$ is unique and, consequently, $F$ is unique and its
 distribution
 inverse is given by (\ref{10}). Finally, observe that
 with this choice of $F^{-1}$, the equality is also
 attained in Moriguti's inequality, because
 $F^{-1}$ is constant in the interval where
 $\overline{G}_{k:n}<G_{k:n}$.
 $\Box$
 \bigskip
 \end{pr}

 \begin{REM}
 \label{rem.4}
 For $k=n\geq 2$, $G_{n:n}(u)=\overline{G}_{n:n}(u)=u^n$ and
 $g_{n:n}(u)=\overline{g}_{n:n}(u)=n u^{n-1}$. Therefore,
 the optimal population is given by
 \[
 F^{-1}(u)
 =\mu \Big(\frac{n-\alpha}{1-\alpha}\Big) u^\frac{n-1}{1-\alpha},
 \ \ 0<u<1,
 \]
 and this corresponds to a power-type distribution function:
 \[
 F(x)=\left(\frac{(1-\alpha)x}{(n-\alpha)\mu}\right)^\frac{1-\alpha}{n-1}, \ \
 0\leq x\leq \frac{n-\alpha}{1-\alpha} \mu.
 \]
 Moreover, the optimal bound for the maximum, (\ref{8}),
 reads as
 \[
 \E (X_{n:n})^\alpha \leq n\Big(\frac{1-\alpha}{n-\alpha}\Big)^{1-\alpha} \
 \mu^\alpha.
 \]
 It is worth pointing out that $\lim_{\alpha\nearrow 1}
 n\Big(\frac{1-\alpha}{n-\alpha}\Big)^{1-\alpha}=n$, yielding
 the best possible non-attainable
 bound $\E X_{n:n}\leq n\mu$; see Corollary \ref{cor.1}.
 \end{REM}

 \begin{REM}
 Due to a result of Balakrishnan (1993),
 the value of $\rho_{2:n}$
 can be calculated in a closed form. In fact,
 $\rho_{2:n}=1/(n-1)^2$ and, consequently,
 \[
 \overline{g}_{2:n}(u)=
 \left\{
 \begin{array}{ll}
 n(n-1) u (1-u)^{n-2}, & 0<u\leq \frac{1}{(n-1)^2},
 \\
 n(n-1)\rho_{2:n}(1-\rho_{2:n})^{n-2}=\frac{n^{n-1}(n-2)^{n-2}}{(n-1)^{2n-3}}, &
 \frac{1}{(n-1)^2}\leq u<1.
 \end{array}
 \right.
 \]
 Hence, for $n\geq 3$, (\ref{9}) reads as
 \[
 A_{2:n}(\alpha)=n(n-1)\left\{\rho_{2:n}^\frac{1}{1-\alpha} (1-\rho_{2:n})^
 \frac{n-1-\alpha}{1-\alpha}+\int_{0}^{\rho_{2:n}}
 u^\frac{1}{1-\alpha}(1-u)^\frac{n-2}{1-\alpha} du
 \right\}^{1-\alpha}, \ \ 0<\alpha<1.
 \]
 This expression should be compared to the corresponding one
 in Example \ref{exam.1}, highlighting the different nature
 of the cases $\alpha\leq 1$, $\alpha\geq 1$.
 \end{REM}

 \section{The minimum in the independent and in the iid case}
 So far, we have not examine the minimum. The reason is that for
 the minimum we can obtain somewhat more complete results,
 that is, sharp upper bounds, including the independent case.
 We start with the iid case.

 \begin{theo}
 \label{th1}
 Let $X$ be a nonnegative rv with
 $\E X=\mu\in(0,\infty)$, and assume that $X_1,\ldots,X_n$ {\rm ($n\geq 2$)}
 are iid rv's distributed like
 $X$. Then,
 the random variable $X_{1:n}=\min\{X_1,\ldots,X_n\}$
 has finite $n$-th moment. Moreover, the inequality
 \be
 \label{1.1}
 \E (X_{1:n})^n \leq \mu^n
 \ee
 holds true, and the equality is attained if and only if
 there exists a number $p\in(0,1]$ such that $X$
 assumes the values $0$ and $\mu/p$ with respective
 probabilities $1-p$ and $p$.
 \end{theo}
 \begin{pr}{Proof}
 See Theorem \ref{th2}, below, for a more general result.
 $\Box$
 \end{pr}

 \begin{REM}
 \label{rem1}
 Theorem \ref{th1} can be viewed in another form, as follows:
 If $X$  is a nonnegative rv with
 $\E X^{1/n}<\infty$ for some $n$, then
 the minimum $X_{1:N}$ is integrable for all $N\geq n$, and, moreover,
 \[
  \E X_{1:N} \leq \E X_{1:n} \leq \left(\E X^{1/n} \right)^n, \ \ N\geq n.
 \]
 Note that, for any  $N\geq n$,
 the upper bound $\left(\E X^{1/n}\right)^n$
 is best possible for $\E X_{1:N}$; this happens
 because we did not exclude a degenerate rv $X$.
 \end{REM}
 \begin{REM}
 \label{rem2}
 The result of Theorem \ref{th1} cannot be extended to any
 higher moment (except in the trivial case $\mu=0$); see
 Lemma \ref{lem.2}.
 A somewhat more direct computation is as follows:
 consider the rv $X$ with
 df
 \[
 F(x)=1-\frac{e}{x(\log x)^2}, \ \ x\geq e.
 \]
 Using the well-known
 formula
 \be
 \label{1.2}
 \E X = \int_{0}^{\infty} \big(1-F(t)\big) dt,
 \ee
 which is valid for any nonnegative rv,
 it is easily seen that $\E X =2e<\infty$. Also,
 since for any
 $\delta\in(0,\infty)$ and $n\in\{2,3,\ldots\}$,
 the df of $(X_{1:n})^{n+\delta}$ is
 $1-\left(1-F\big(t^{1/(n+\delta)}\big)\right)^n$,
 $t\geq 0$, (\ref{1.2}) yields
 \[
 \E (X_{1:n})^{n+\delta}
 =
 (n+\delta) \int_0^{\infty} x^{n+\delta-1} \big(1-F(x)\big)^n dx
 \geq \int_{e}^{\infty} \frac{(n+\delta)e^n}{x^{1-\delta}(\log x)^{2n}} dx
 =\infty.
 \]
 It is clear that for arbitrary $\mu>0$, the rv $Y=\mu X/(2e)\geq 0$
 has mean $\mu$, and the rv $(Y_{1:n})^{n+\delta}$ is non-integrable for any
 $\delta\in(0,\infty)$ and for any $n\in\{2,3,\ldots\}$.
 \end{REM}
 \begin{REM}
 \label{rem3}
 Theorem \ref{th1} yields the best upper bound for any fractional moment
 of $X_{1:n}$ as follows:
 Since $x\mapsto x^p$ ($0<p<1$) is concave in $[0,\infty)$, Jensen
 (or Lyapounov) inequality, combined with (\ref{1.1}),
 yields
 \[
 \E (X_{1:n})^\alpha = \E \Big[ ((X_{1:n})^n)^{\alpha/n}\Big]
 \leq \Big(\E (X_{1:n})^n \Big)^{\alpha/n}\leq \mu^\alpha,  \ \ 0<\alpha\leq n.
 \]
 The upper bound $\mu^\alpha$ is clearly best possible, since it is
 attained (uniquely,
 unless $\alpha=n$)
 by a degenerate $X$ at $\mu$.
 \end{REM}

 From now on, assume that $X_1,\ldots,X_n$ are independent, non-negative
 rv's
 with finite means $\E X_i=\mu_i>0$ for $i=1,\ldots,n$, and set
 $X_{1:n}=\min\{X_1,\ldots,X_n\}$. Our purpose is to derive
 the best possible upper bounds for the moments of
 $X_{1:n}=\min\{X_1,\ldots,X_n\}$, and
 the populations that attain the bounds.
 \begin{theo}
 \label{th2}
 {\rm (i)} The random variable $X_{1:n}$ has finite $n$-th moment and,
 moreover,
 the inequality
 \be
 \label{2.1}
 \E X_{1:n}^n \leq \mu_1\cdots \mu_n
 \ee
 is valid, with equality if and only if there exists a number
 $M\geq \max_i\{\mu_i\}$
 such that
 \[
 \Pr(X_i=M)=\frac{\mu_i}{M}=1-\Pr(X_i=0), \ \ \ i=1,\ldots,n.
 \]
 {\rm (ii)} For given strictly positive numbers
 $\mu_1,\ldots,\mu_n$ {\rm ($n\geq 2$)} we can find independent,
 non-negative r.v.'s $X_1,\ldots,X_n$ such that
 \[
 \E X_i=\mu_i \ (i=1,\ldots,n) \ \ \mbox{and }  \ \
 \E (X_{1:n})^{n+\delta} =\infty \ \ \mbox{for all } \delta\in(0,\infty).
 \]
 Furthermore, if $\mu_1=\cdots=\mu_n$, the r.v.'s $X_1,\ldots,X_n$ can
 be chosen to be iid.
 \end{theo}
 \begin{pr}{Proof}
 (i)
 To see the inequality (\ref{2.1}), just take expectations
 to the obvious (deterministic) inequality
 \[
 (X_{1:n})^n
 \leq X_1\cdots X_n.
 \]
 Moreover, observe that $X_1\cdots X_n=X_{1:n}\cdots X_{n:n}$,
 where $X_{1:n}\leq \cdots \leq X_{n:n}$ are the corresponding order
 statistics of $X_1,\ldots,X_n$.
 Thus,
 for the equality to hold, it is necessary and sufficient
 that
 \[
 \E\Big[X_{1:n}\Big(X_{2:n}\cdots X_{n:n}-(X_{1:n})^{n-1}\Big)\Big]=0.
 \]
 This implies the relation
 \be
 \label{2.2}
 \Pr\Big(\{X_{1:n}=0\} \cup \{X_1=\cdots=X_n>0\}\Big) =1.
 \ee
 Let $p_i=\Pr(X_i>0)>0$ (since $\mu_i>0$). It follows that
 $\Pr(X_{1:n}>0)=\prod_{i=1}^n p_i>0$ and, from (\ref{2.2}),
 \be
 \label{2.3}
  \Pr(X_1=\cdots=X_n>0 \  |  \ X_{1:n}>0)=1.
 \ee
 Define now the independent rv's $Y_i$ with $Y_i\law (X_i\ | \ X_i>0)$;
 that is, $F_{Y_i}(y)=(F_{X_i}(y)-1+p_i)/p_i$, $y\geq 0$. Then
 (\ref{2.3}) reads as $\Pr(Y_1=\cdots=Y_n)=1$ and,
 by the independence of $Y_i$,  it follows
 that we can find a constant $M>0$ such that $\Pr(Y_i=M)=1$ for
 all $i$; hence, $\Pr(X_i=0)+\Pr(X_i=M)=1$. From $\mu_i=\E X_i=M p_i$
 we get $p_i=\mu_i/M$ and, thus, $M\geq \max_i\{\mu_i\}$.
 As a final check, it is easily verified that the rv's $X_i$
 with $\Pr(X_i=M)=\mu_i/M=1-\Pr(X_i=0)$ attain the equality in
 (\ref{2.1}).

 For the proof of (ii), see Lemma \ref{lem.2}.
 \bigskip
 $\Box$
 \end{pr}

 It became clear from Theorem \ref{th2}(ii) and Remark \ref{rem2}
 that we cannot hope for finiteness of moments of order higher than $n$
 (for $X_{1:n}$) without additional assumptions. It is, thus, desirable,
 to derive upper bounds for lower moments. Indeed, in this case we
 have the following result.
 \begin{theo}
 \label{th3}
 Let $X_1,\ldots,X_n$ be independent, non-negative,
 rv's with finite expectations $\E X_i=\mu_i>0$ and, without
 loss of generality, assume that
 $0<\mu_1\leq \mu_2\leq \cdots\leq \mu_n$.
 Then, for every $\alpha\in(0,n]$ we have
 \be
 \label{2.7}
 \E (X_{1:n})^\alpha
 \leq \mu_1\cdots \mu_{k-1} (\mu_k)^{\alpha-k+1}, \ \ \alpha\in(k-1,k],
 \ k=1,\ldots,n.
 \ee
 The bound is best possible, since the equality is attained
 by the independent rv's $X_i$ with
 \be
 \label{2.8}
 \begin{array}{ll}
 \Pr(X_i=\mu_k)=\frac{\mu_i}{\mu_k}=1-\Pr(X_i=0),
 & i=1,\ldots, k,
 \\
 \Pr(X_i=\mu_i)=1, &
 i=k+1,\ldots,n,
 \end{array}
 \ee
 where $k\in\{1,\ldots,n\}$ is the unique integer
 such that $k-1<\alpha\leq k$.
 \end{theo}
 \begin{pr}{Proof}
 Since it is easily checked that the rv's in (\ref{2.8})
 attain the equality in (\ref{2.7}), we proceed to
 verify the inequality (\ref{2.7}). To this end,
 fix $\alpha\in(k-1,k]$
 and consider
 the following deterministic inequalities,
 valid for $X_i\geq 0$:
 \begin{eqnarray*}
 \min\{X_1,\ldots,X_n\}
 &\leq & X_1
 \\
 \min\{X_1,\ldots,X_n\}
 &\leq & X_2
 \\
 &\vdots&
 \\
 \min\{X_1,\ldots,X_n\}
 &\leq & X_{k-1}
 \\
 (\min\{X_1,\ldots,X_n\})^{\alpha-(k-1)}
 &\leq & (X_k)^{\alpha-(k-1)}.
 \end{eqnarray*}
 Multiplying, we get
 \be
 \label{2.9}
 (X_{1:n})^\alpha\leq X_1\cdots X_{k-1}(X_k)^{\alpha-k+1}.
 \ee
 Hence, taking expectations in (\ref{2.9})
 and using independence,
 we deduce the inequality
 \[
 \E (X_{1:n})^\alpha \leq \mu_1\cdots\mu_{k-1}\E (X_{k})^{\alpha-k+1}.
 \]
 Finally, since $0<\alpha-k+1\leq 1$,
 the function $x\mapsto x^{\alpha-k+1}$ is concave in $[0,\infty)$,
 and Jensen (or Lyapounov) inequality yields (\ref{2.7}).
 $\Box$
 \bigskip
 \end{pr}

 Notice that the inequality (\ref{2.9}) shows that $(X_{1:n})^\alpha$
 (for $\alpha\in(k-1,k]$) is integrable even if $\mu_{k+1}=\infty$;
 this is explained
 from the fact that $X_{1:n}\leq \min\{X_1,\ldots,X_k\}$
 and, by Theorem \ref{th2},
 $X_{1:k}$ has finite $k$-th (hence $\alpha$-th) moment.
 Note also that (\ref{2.7}) yields Remark \ref{rem3}
 for the iid case.
 \bigskip

 \noindent
 {\bf Acknowledgements.} I would like to thank Dimitris Cheliotis
 for proposing the problem, and for providing me with a
 partial answer for the case $n=2$.

 {
 
 }

 \end{document}